\documentclass[a4paper,12pt]{amsart}
\usepackage{fullpage}
\usepackage[english,francais]{babel}
\usepackage{epsf}
\usepackage[dvips]{epsfig}
\usepackage{amssymb}
\usepackage{amsmath}
\usepackage[all]{xy}
\usepackage{color}

\usepackage{ucs}
\usepackage[utf8x]{inputenc}\usepackage{wasysym}
\usepackage{aeguill}

\newtheorem{theorem}{Theorem}[section]
\newtheorem{prop-f}[theoreme]{Proposition}
\newtheorem{prop}[theorem]{Proposition}

\newtheorem{corollary}[theorem]{Corollary}

\newtheorem{lemma}[theorem]{Lemma}

\newcommand{\finpreuve}{\hspace{\stretch{1}}{$\square$}}

\newcommand{\R}{\mathbb{R}}

\newcommand{\Z}{\mathbb{Z}}
\newcommand{\N}{\mathbb{N}}

\def\esp{\medskip\noindent}
\def\prg#1{\esp{\bf #1. }}

\def\proof{\prg{Proof}}

\def\proofof#1{\prg{Proof of #1}}
\renewcommand{\epsilon}{\varepsilon}
\def\btab{\begin{eqnarray*}}
\def\etab{\end{eqnarray*}}
\def\beq{\begin{equation}}
\def\eeq{\end{equation}}

\newcounter{numeroexo}

\begin{document}

\selectlanguage{english}

%\title[Continuum percolation in high dimensions]{Continuum percolation in high dimensions}
\title[Continuum percolation in high dimensions]{Non-optimality of constant radii in high dimensional continuum percolation}
{
\author{Jean-Baptiste Gou\'er\'e}
\author{R{\'e}gine Marchand}
\address{Laboratoire de Math{\'e}matiques, Applications et Physique
Math{\'e}matique d'Orl{\'e}ans UMR 6628\\ Universit{\'e} d'Orl{\'e}ans\\ B.P.
6759\\
 45067 Orl{\'e}ans Cedex 2 France}
\email{Jean-Baptiste.Gouere@univ-orleans.fr}

\address{1. Universit{\'e} de Lorraine\\
Institut Elie Cartan de Lorraine, UMR 7502 (math{\'e}matiques)\\
Vandoeuvre-l{\`e}s-Nancy, F-54506 France. 
2. CNRS \\
Institut Elie Cartan de Lorraine, UMR 7502 (math{\'e}matiques)\\
Vandoeuvre-l{\`e}s-Nancy, F-54506 France. }
\email{Regine.Marchand@univ-lorraine.fr}
}

%\subjclass[2000]{60K35, 82B43.} 
%\keywords{\motsclefs}
% ici pour article 

\begin{abstract}
Consider a Boolean model $\Sigma$ in $\R^d$.
The centers are given by a homogeneous Poisson point process with intensity $\lambda$ and the radii of distinct balls are i.i.d.\ with common distribution $\nu$.
The critical covered volume is the proportion of space covered by $\Sigma$ when the intensity $\lambda$ is critical for percolation.
Previous numerical simulations and heuristic arguments suggest that the critical covered volume
may be minimal when $\nu$ is a Dirac measure. 
In this paper, we prove that it is not the case in sufficiently high dimension.
\end{abstract}

\maketitle

\section{Introduction and statement of the main results}

The Boolean model is a popular model for continuum percolation. It can be described in the following way. Let $\nu$ be a finite measure on $(0,+\infty)$, with positive mass.
Let $d \ge 2$ be an integer, $\lambda > 0$ be a real number and  
$\xi$ be a Poisson point process on $\R^d \times (0,+\infty)$ 
whose intensity measure is the Lebesgue measure on $\R^d$ times $\lambda\nu$.
The Boolean model $\Sigma(\lambda\nu)$ in $\R^d$ driven by $\lambda\nu$ is the following random subset of  $\R^d$:
$$
\Sigma(\lambda\nu) = \bigcup_{(c,r) \in \xi} B(c,r),
$$
where $B(c,r)$ is the open Euclidean ball centered at $c \in \R^d$ and with radius $r \in (0,+\infty)$. Note that the collection of centers of the balls of the Boolean model is a homogeneous Poisson point process on $\R^d$ with intensity $\lambda \nu((0,+\infty))$, and that the radii of the distinct balls are i.i.d.\ with law $\nu(.)/\nu((0,+\infty))$, 
and independent of the point process of the centers. In our study, we  focus on the Boolean model with deterministic radii (when $\nu$ is a Dirac mass $\delta_\rho$, with $\rho>0$) and on the Boolean model with two distinct radii (when $\nu$ is a weighted sum of two Dirac masses).

We say that $\Sigma(\lambda\nu)$ percolates if the probability that there is an unbounded connected component of $\Sigma(\lambda\nu)$ that contains the origin is positive.
This is equivalent to the almost-sure existence of an unbounded connected component of $\Sigma(\lambda\nu)$.
We refer to the book by Meester and Roy \cite{Meester-Roy-livre} for background on continuum percolation.
The critical intensity is defined by:
$$
\lambda^c_d(\nu) = \inf \{\lambda>0 : \Sigma(\lambda\nu) \hbox{ percolates}\}.
$$
%One easily checks the following:
%\begin{itemize}
%\item If $\lambda<\lambda^c(\nu)$ then, with probability one, all the connected components of $\Sigma(\lambda\nu)$ are bounded.
%\item If $\lambda>\lambda^c(\nu)$ then, with probability one, there exists one connected component in $\Sigmma(\lambda\nu)$.
%\end{itemize}
One easily checks that $\lambda^c_d(\nu)$ is finite, and in
 \cite{G-perco-boolean-model} it is proven that $\lambda^c_d(\nu)$ is positive if and only~if 
\begin{equation}\label{e:volumefini}
\int r^d\nu(dr) < +\infty. 
\end{equation}
We assume that this assumption is fulfilled.

By ergodicity, the Boolean model $\Sigma(\lambda\nu)$ has a deterministic natural density.
This is also the probability that a given point belongs to the Boolean model and  it is given by :
$$
P(0 \in \Sigma(\lambda\nu))=1-\exp\left(-\lambda \int v_dr^d \nu(dr)\right),
$$
where $v_d$ denotes the volume of the unit ball in $\R^d$.
The critical covered volume $c_d^c(\nu)$ is the density of the Boolean model when the intensity is critical :
$$
c^c_d(\nu) =1-\exp\left(-\lambda_d^c(\nu) \int v_dr^d \nu(dr) \right).
$$
Unlike the critical intensity $\lambda^c_d$, the critical covered volume $c_d^c$ is invariant under scaling. For all $a>0$, let $H^a(\nu)$ be the image of $\nu$ under the map defined by $x \mapsto ax$.
We have the following scaling property:
\begin{equation}\label{e:scaling}
c^c_d(H^a\nu)=c^c_d(\nu).
\end{equation}
Indeed, a critical Boolean model remains critical when rescaling and the density is invariant by rescaling \footnote{By rescaling, we mean multiplying all coordinates and radii by the same scalar.}. More formally, this invariance of the critical covered volume $c_d^c$ under rescaling is a consequence of Proposition 2.11 in \cite{Meester-Roy-livre}.
Note for example that for any $\rho>0$, 
$$c_d^c(\delta_1)=c_d^c(\delta_\rho), \text{ while } \lambda_d^c(\delta_1)=\rho^d \lambda_d^c(\delta_\rho).$$ 
One also easily checks the following invariance property: for all $a>0$, $c_d^c(a\nu)=c_d^c(\nu)$.

Practically, we study the critical covered volume through the normalized critical intensity:
$$
\widetilde{\lambda}^c_d(\nu) = \lambda_d^c(\nu) \int v_d(2r)^d \nu(dr).
$$
We then have ${c^c_d}(\nu) =1-\exp\left(-\frac{\widetilde\lambda_d^c(\nu)}{2^d} \right)$.
The factor $2^d$ may seem arbitrary here, 
its interest will appear in the statement of the next theorems. Note also that the normalized critical intensity $\widetilde{\lambda}^c_d$ is also invariant under rescaling. 

\subsection*{Normalized critical intensity as a function of $\nu$.}
It has been conjectured by Kert\'esz and Vicsek~\cite{Kertesz-al} that the normalized critical intensity should be independent of $\nu$,
as soon as the support of $\nu$ is bounded.
Phani and Dhar \cite{Phani-Dhar} gave a heuristic argument suggesting that the conjecture were false.
A rigorous proof was then given by Meester, Roy and Sarkar in \cite{Meester-Roy-Sarkar}.
More precisely, they gave examples of measures $\nu$ with two atoms such that:
\begin{equation}\label{e:MRS}
\widetilde\lambda^c_d(\nu) > \widetilde\lambda_d^c(\delta_1).
\end{equation}
As a consequence of Theorem 1.1 in the paper by Menshikov, Popov and Vachkovskaia~\cite{Menshikov-al-multi},
we even get that $\widetilde\lambda^c_d(\nu)$ can be arbitrarily large\footnote{Actually the result of \cite{Menshikov-al-multi} is a much stronger statement than the consequence we use here.}.
%More precisely, if
%\begin{equation}\label{e:cvmulti}
%\text{if }\nu(n,a)=\sum_{k=0}^{n-1} a^{dk} \delta_{a^{-k}}, \text{ then }\widetilde\lambda^c_d(\nu(n,a)) \to n\widetilde\lambda^c_d(\delta_1) \hbox{ as } a \to \infty.
%\end{equation}
%Actually the result of \cite{Menshikov-al-multi} is the following much stronger statement: 
%$\lambda^c_d(\nu(+\infty,a)) \to \lambda_d^c(\delta_1)$ when $a \to \infty$.
%The convergence \eqref{e:cvmulti} is implicit in the work of Meester, Roy and Sarkar in \cite{Meester-Roy-Sarkar}, 
%at least when $n=2$.
%There were also heuristics for such a result in \cite{Phani-Dhar}.
On the contrary,  Theorem 2.1 in \cite{G-perco-boolean-model} gives the existence of a positive constant $C_d$, that depends only on the 
dimension $d$, such that, for all $\nu$ satisfying \eqref{e:volumefini}:
$$
\widetilde\lambda^c_d(\nu) \ge C_d.
$$
To sum up, $\widetilde\lambda_d^c(\cdot)$ is not bounded from above but is bounded from below by a positive constant.
In other words, the critical covered volume $c^c_d(\cdot) \in (0,1)$ can be arbitrarily close to $1$ 
but is bounded from below by a positive constant. 
It is thus natural to seek optimal measures, that is the ones which minimize the normalized critical intensity, or equivalently, the critical covered volume.

\begin{figure}[h!]
\includegraphics[scale=0.7]{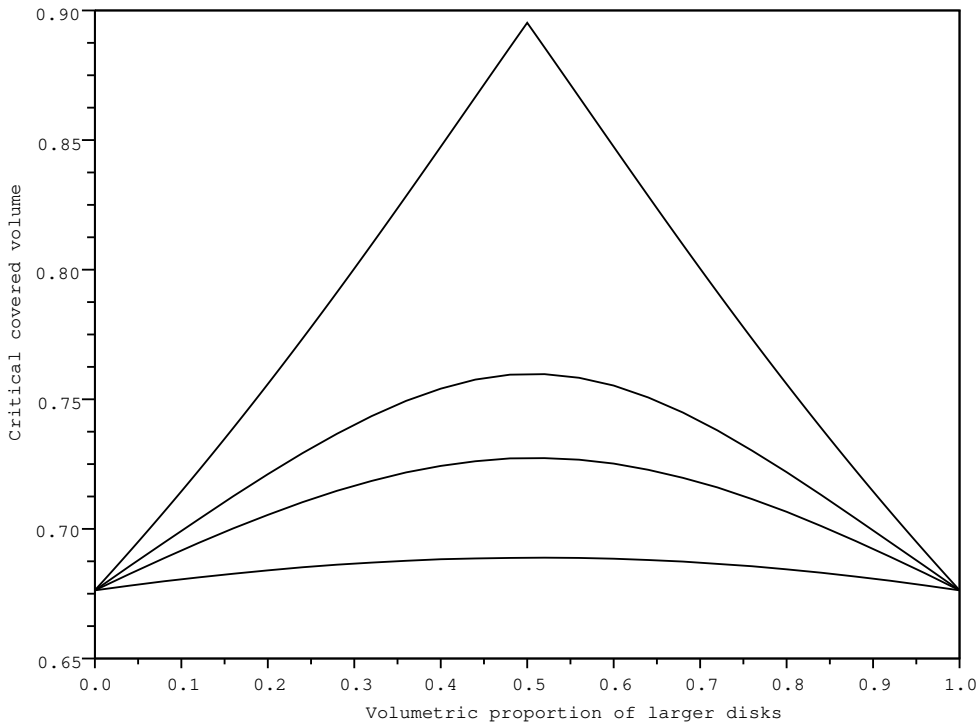}
\caption{Critical covered volume as a function of $\alpha$ for different values of $\rho$.
>From bottom to top: $\rho=2, \rho=5, \rho=10$ and the limit as $\rho \to \infty$.}
\label{f:multiscale}
\end{figure}

In the physical literature, it is strongly believed that, at least when $d=2$ and $d=3$, the critical covered volume is minimum in the case of a deterministic radius, that is when the distribution of radii is a Dirac measure.
This conjecture is supported by numerical evidence (to the best of our knowledge, the most accurate estimations are given in
a paper by Quintanilla and Ziff \cite{QZ-PRE-2007} when $d=2$ and in
a paper by Consiglio, Baker, Paul and Stanley \cite{Consiglio-2003} when $d=3$). On Figure \ref{f:multiscale}, we plot the critical covered volume in dimension 2 as a function of $\alpha$ and for different values of $\rho$ when $\nu=(1-\alpha) \delta_1+{\alpha}{\rho^{-2}}\delta_\rho$. The data for finite values of $\rho$ come from numerical estimations in \cite{QZ-PRE-2007}, while the data for the limit of $\rho$ going to infinity come from the study of the multi-scale Boolean model in \cite{G-multi}. See Section 1.4 in \cite{G-multi} for further references. 
The conjecture is also supported by some heuristic arguments in any dimension (see for example Dhar~\cite{Dhar-1997}, and Balram and Dhar~\cite{Balram-Dhar}).
In \cite{Meester-Roy-Sarkar}, 
it is noted that the rigorous proof of \eqref{e:MRS} suggests that the deterministic case might be optimal for any $d \ge 2$.
%\verifier{en physique trouver des valeurs approchées au moins dans des régimes asymptotiques a aussi donné lieu \C3  quelques papiers}

In this paper we show on the contrary that for all $d$ large enough the critical covered volume is not minimized by the case of deterministic radii. 

\subsection*{Normalized critical intensity  in high dimension : the case of a deterministic radius.}

Assume here that the measure $\nu$ is a Dirac mass at $1$, that is that the radii of the balls are all equal to $1$. 
Penrose proved the following result in \cite{Penrose-high-dimensions} :
\begin{theorem}[Penrose] \label{t:penrose}
$\displaystyle\lim_{d\to\infty} \widetilde\lambda^c_d(\delta_1) =1.$
\end{theorem}
With the scale invariance of $\widetilde\lambda^c_d$, this limit can readily be generalized to any constant radius : for any $\rho>0$,
$$
%\lim_{d\to\infty} \widetilde\lambda^c_d((\alpha \delta_\rho)_d)=
\lim_{d\to\infty} \widetilde\lambda^c_d(\delta_\rho) = 
\lim_{d\to\infty} \widetilde\lambda^c_d(\delta_1) =1.
$$   
Theorem \ref{t:penrose} is the continuum analogue of a result of Kesten
\cite{Kesten-high-dimensions} 
for Bernoulli bond percolation on the nearest-neighbor integer lattice $\Z^d$, which says that the critical percolation parameter is asymptotically equivalent to $1/(2d)$.

Let us say a word about the ideas of the proof of Theorem \ref{t:penrose}.

The inequality $\widetilde\lambda^c_d(\delta_1) > 1$  
holds  for any $d \ge 2$.
The proof is simple, and here is the idea. We consider the following natural genealogy.
The deterministic ball $B(0,1)$ is said to be the ball of generation $0$.
The random balls of $\Sigma(\lambda\delta_1)$ that touch $B(0,1)$ are then the balls of generation $1$.
The random balls that touch one ball of generation $1$ without being one of them are then  the balls of generation $2$ and so on.
Let us denote by $N_d$ the number of all balls that are descendants of $B(0,1)$.
There is no percolation if and only if $N_d$ is almost surely finite.

Now denote by $m$ the Poisson distribution with mean $\lambda v_d 2^d $ : this is the law of the
number of balls of $\Sigma(\lambda \delta_1)$ that touch a given ball of radius $1$.
Therefore, if there were no interference between children of different balls, $N_d$ would be equal to $Z$, 
the total population in a Galton-Watson process with offspring distribution $m$.
Because of the interferences due to the fact that the Boolean model lives in $\R^d$, this is not true : in fact, $N_d$ is only stochastically dominated by $Z$.
Therefore, if $\lambda v_d2^d \le 1$, then $Z$ is finite almost surely, so $N_d$ is finite almost surely
and therefore there is no percolation.
This implies 
$$\widetilde\lambda_d^c(\delta_1) = v_d2^d \lambda_d^c(\delta_1) > 1.$$
The difficult part of Theorem \ref{t:penrose} is to prove that if $d$ is large, then the interferences are small,
so $N_d$ is close to $Z$ and therefore there is percolation for large $d$ as soon as 
$v_d2^d \lambda$ is a constant strictly larger than one.

To sum up, at first order, the asymptotic behavior of the critical intensity of the Boolean model with constant radius is given by 
the threshold of the associated Galton-Watson process, as in the case of Bernoulli percolation on $\Z^d$ : roughly speaking, as the dimension increases, the geometrical constraints of the finite dimension space decrease and at the limit, we recover the non-geometrical case of the corresponding Galton-Watson process.

\subsection*{Normalized critical intensity  in high dimension : the case of radii taking two values.}

Let $1 < \rho < 2$.
Set 
$$
\mu = \delta_1+ \delta_\rho.
$$
If $d\ge 1$ is an integer, we define the normalized measure $\mu_d$ on $(0,+\infty)$ by setting :
\begin{equation}
\label{e:mud}
\mu_d= \delta_1+\frac{1}{\rho^d} \delta_\rho.
\end{equation}
We will study the behavior of $\widetilde\lambda_d^c(\mu_d)$ as $d$ tends to infinity.
Let us motivate the definition of $\mu_d$ with the following two related properties:
\begin{enumerate}
\item
Consider the Boolean model $\Sigma(\lambda\mu_d)$ on $\R^d$ driven by $\lambda\mu_d$ where $\lambda>0$.
The number of balls of $\Sigma(\lambda\mu_d)$ with radius $1$ that contains a given point
is a Poisson random variable with intensity $\lambda v_d$. 
The number of balls of $\Sigma(\lambda\mu_d)$ with radius $\rho$ that contains a given point
is also a Poisson random variable with intensity $\lambda v_d$. 
Loosely speaking, this means that contrary to what happens in the Boolean model driven by $\lambda \mu$, the relative importance of the two types of radii does not depend on the dimension $d$ in the Boolean model driven by $\lambda \mu_d$.
\item
A closely related property is the following one.
Consider two independent Boolean model $\Sigma$ and $\Sigma'$,
both driven by $\lambda \delta_1$.
Then $\Sigma \cup \rho \Sigma'$ is a Boolean model driven by $\lambda\mu_d$.
\end{enumerate}

\begin{theorem} \label{t:2}
Let $1<\rho<2$.
Set as before $\mu_d = \delta_1+ \rho^{-d}\delta_\rho$.
Then
\begin{equation}\label{e:2}
\lim_{d \to +\infty} \frac1d \ln\left(\widetilde{\lambda}^c_d(\mu_d)\right) = \ln\left(\kappa_\rho^c\right), \text{ where } \kappa_\rho^c=\frac{2 \sqrt \rho}{1+\rho}. \nonumber
\end{equation}
\end{theorem}	

Note that as $1<\rho<2$, $\kappa_\rho^c<1$. The following result is then an immediate consequence of Theorem \ref{t:penrose} and Theorem \ref{t:2}.

\begin{corollary} \label{t:conjecture}
If the dimension $d$ is large enough, then there exists a probability measure $\nu$ on $(0,+\infty)$ such that :
$$
c_d^c(\nu) < c_d^c(\delta_1).
$$
In other words, the conjecture is false in high dimensions.
\end{corollary}

We end this section by some remarks:
\begin{itemize}
\item One can easily extend Theorem \ref{t:2} as follows. 
Let $\alpha, \beta, a, b >0$. 
Set $\rho=b/a$ and assume $1<\rho<2$.
Then
$$
\lim_{d \to +\infty} \frac1d \ln\left(\widetilde{\lambda}^c_d(\alpha a^{-d} \delta_a + \beta b^{-d} \delta_b)\right) = \ln\left(\kappa_\rho^c\right)<0.
$$
\item As we will see in the proof, the critical threshold $\kappa_\rho^c$ is given by the critical parameter of an associated two-types Galton-Watson process when $1 < \rho <2$; we  prove in a companion paper \cite{GouMar-Bool-Perco-high-dim} that this is not the case for $\rho>2$. 
\item If one does not normalize the distribution one has 
\footnote{The upper bound can be proven using $\lambda_d^c(\alpha \delta_a+\beta \delta_b) \le \lambda_d^c(\beta \delta_b)$.
The lower bound can be proven using the easy part of the comparison with a two-type Galton-Watson process.}
$\widetilde{\lambda}_d^c(\alpha \delta_a+\beta \delta_b) \to 1$ and thus 
$\widetilde{\lambda}_d^c(\alpha \delta_a+\beta \delta_b) \sim \widetilde{\lambda}_d^c(\delta_1)$.
This behavior is due to the fact that, without normalization, the influence of the small balls vanishes in high dimension. 
\end{itemize}

\section{Proofs}

\subsection{Notations}
\label{s:notations}
Fix $1<\rho<2$ and $\kappa>0$.
Once the dimension $d \ge 1$ is given, we consider two independent stationary Poisson point processes on $\R^d$: $\chi_1$ and $\chi_\rho$, with respective intensities 
$$
\lambda_1=\frac{\kappa^d}{v_d 2^d} \quad \text{ and } \quad 
\lambda_{\rho}=\frac{\kappa^d}{v_d 2^d \rho^d}.
$$
%Thus, $\mu_d= \lambda_1\delta_1+\lambda_{\rho}\delta_\rho$.
To $\chi_1$ and $\chi_{\rho}$, we respectively associate the two Boolean models
$$
\Sigma_1 = \bigcup_{x \in \chi_1} B(x,1) \quad \text{ and } \quad 
\Sigma_{\rho} = \bigcup_{x \in \chi_{\rho}} B(x,\rho).
$$
We focus on the percolation properties of the following two-type Boolean model
$$
\Sigma = \Sigma_1 \cup \Sigma_{\rho}.
$$
This Boolean model is driven by the measure 
$$
\lambda_1 \delta_1 + \lambda_\rho \delta_\rho = \frac{\kappa^d}{v_d 2^d} \mu_d
$$
where $\mu_d$ is defined as before by \eqref{e:mud}. Remember that
$$\kappa_\rho^c=\frac{2 \sqrt \rho}{1+\rho}<1.$$

\subsection{Subcritical phase}

The aim of this subsection is to prove the following result.

\begin{prop} \label{p:subcritical} 
%Recall $1 < \rho < 2$.
If $\kappa<\kappa^c_\rho$, then, as soon as the dimension $d$ is large enough, percolation does not occur in the two-type Boolean model $\Sigma$.
\end{prop}

\proof The proof is very similar to the easy part of the proof of Theorem \ref{t:penrose}.
The only difference is that we consider a two-types Galton-Watson process instead of a one-type Galton-Watson process.
Therefore, we only sketch the proof and refer to \cite{Penrose-high-dimensions} for a more detailed proof.

The idea is to consider the following natural genealogy. 
The deterministic ball $B(0,\rho)$ is said to be the ball of generation $0$.
The random balls of $\Sigma$ that touch $B(0,\rho)$ are then the balls of generation $1$.
They can be of two different types: either of radius $1$ or of radius $\rho$.
The random balls that touch one ball of generation $1$ without being one of them are then  the balls of generation $2$ and so on. 
 
This genealogical process is stochastically dominated by a two-types Gatson-Watson process.
Basically, the Galton-Watson process is obtained by neglecting 
the geometrical constraints due to the fact that the Boolean model lives in $\mathbb R^d$.
It is defined as follows.
Start with one individual of type $\rho$.
The offspring distribution of type $1$ of an individual of type $\rho$ is defined to be the distribution of the number of balls of $\Sigma_1$
that intersect a given deterministic ball of radius $\rho$.
Therefore, it is a Poisson random variable with mean $\lambda_1 v_d(1+\rho)^d$.
The other offspring distributions are defined similarly. 
The matrix of means of offspring distributions is thus given by:
$$
M_d=
\begin{pmatrix}
\lambda_1 v_d(1+1)^d &  \lambda_\rho  v_d(1+\rho)^d \\
\lambda_1 v_d(1+\rho)^d & \lambda_\rho v_d(\rho+\rho)^d 
\end{pmatrix}= \kappa^d
\begin{pmatrix}
1 &  \left(\frac{1+\rho}{2\rho}\right)^d \\
\left(\frac{1+\rho}{2}\right)^d & 1 
\end{pmatrix} .
$$
Let $r_d$ denote the largest eigenvalue of $M_d$.
The extinction probability of the two-types Galton-Watson process is $1$ if and only if $r_d \le 1$.
We have:
$$
r_d \sim \left(\frac{\kappa(1+\rho)}{2\sqrt{\rho}}\right)^d.
%, \text{ and thus }\kappa = \frac{2\sqrt{\rho}}{1+\rho} \text{ is the critical parameter.}
$$
As $\kappa < \kappa_\rho^c$, we get that the Galton-Watson process is subcritical for large enough $d$.
Therefore, for large enough $d$, the total progeny of the Galton-Watson process is almost surely finite.
Thus, almost surely, there is no infinite cluster of the Boolean model $\Sigma$ that touches  $B(0,\rho)$. 
As a consequence, almost surely, there is no infinite cluster in the Boolean model $\Sigma$.
\finpreuve

%With Theorem \ref{t:2}, we thus see that the comparison with the two-type Galton-Watson is asymptotically 
%sharp on a logarithmic scale when $1 < \rho <2$; we will prove in a companion paper that this is not the case for $\rho>2$. 

\subsection{Supercritical phase}
\label{s:supercritical}

\subsubsection{Result}

For every $n \ge 0$, we set $R_n=\rho$ if $n$ is even and $R_n=1$ otherwise.
We say that alternating percolation  occurs if there exists an infinite sequence  of distinct points $(x_n)_{n \in \N}$ in $\R^d$  such that, for every $n\ge 0$: 
\begin{itemize}
\item $x_n \in \chi_{R_n}$.
\item $B(x_n,R_n) \cap B(x_{n+1},R_{n+1}) \neq  \emptyset$.
\end{itemize}
In other words, alternating percolation occurs if there exists an infinite path along which  balls of radius $1$ alternate with  balls of radius $\rho$.
The aim of this subsection is to prove the following proposition :

\begin{prop} \label{p:supercritical}
%Recall that $1<\rho<2$. 
Assume $\kappa > \kappa_\rho^c$.
If the dimension $d$ is large enough, then alternating percolation  occurs in $\Sigma$ with probability one. 
\end{prop}

By a straightforward coupling argument, one sees that it is sufficient to prove the proposition under the following assumptions on $\kappa$:
$$
\kappa > \kappa^c_{\rho} \hbox{ and } \kappa < \frac{2 \sqrt 2}{1+\rho} \hbox{ and } \kappa < 1.
$$
We make this assumption in the remaining of this subsection.

We will prove that alternating percolation occurs in the two-type Boolean model in the supercritical case  by embedding in the Boolean model a supercritical $2$-dimensional oriented percolation process.

We thus specify the two first coordinates, and introduce the following notations.
When $d \ge 3$, for any $x \in \R^d$, we write 
$$
x=(x',x'') \in \R^2 \times \R^{d-2}.
$$ 
We write $B'(c,r)$ for the open Euclidean ball of $\R^2$ with center $c\in \R^2$ and radius $r>0$.
In the same way we denote by $B''(c,r)$ the open Euclidean ball of $\R^{d-2}$ with center $c \in \R^{d-2}$ and radius $r>0$.

\subsubsection{One step in the $2$-dimensional oriented percolation model}
\label{s:onestep}

The point here is to define the event that will govern the opening of the edges in the $2$-dimensional oriented percolation process : it is naturally linked to the existence of a finite path composed of a ball of radius $1$ and a ball of radius $\rho$.

We define, for a given dimension $d \ge 3$, the two following subsets of $\R^d$:
\begin{eqnarray*}
W & = &  d^{-1/2} \left((-1,1) \times (-1,0) \times \R^{d-2}\right), \\
W^+ & = &  d^{-1/2} \left((0,1) \times (0,1) \times \R^{d-2}\right). 
\end{eqnarray*}
For  $x_0 \in W$ we set :
\begin{equation} \label{e:bonevnt}
{\mathcal G}^+(x_0) = \left\{
\begin{array}{c}
\hbox{There exist } x_1 \in \chi_1 \cap  W^+ \hbox{ and } x_{2} \in \chi_{\rho} \cap W^+ \\
\hbox{such that } B(x_0, \rho)\cap B(x_1,1) \ne \varnothing \text{ and } B(x_1,1)\cap B(x_{2}, \rho) \ne \varnothing
\end{array}
\right\}.
\end{equation}
Our goal here is to prove that the probability of occurrence of this event is asymptotically large :

\begin{prop} \label{p:onestep} 
%Recall $1<\rho<2$. 
Assume that $\kappa \in (\kappa^c_\rho,1)$.
Choose $p \in (0,1)$. 
If the dimension $d$ is large enough, then for every $x_0 \in W$,
$$
P({\mathcal G}^+(x_0))\ge p.
$$
\end{prop}
Note already that by translation invariance, $P({\mathcal G}^+(x_0))$ does not depend on $x''_0$, so we can assume without loss of generality that $x''_0=0$. 
We introduce the following subsets: 
$$
\begin{array}{|c|c|}
\hline
\text{subsets of } \R^2 & \text{subsets of } \R^{d-2} \\
\hline
\hline
D'_0  =  d^{-1/2}(-1,1) \times (-1, 0) 
& C''_0  =  \{0\} \\
\hline
 D'_1= d^{-1/2}(0,1) \times (0,1)
&  C''_1  =  B'' \left(0, (1+\rho)-\frac{6}{d} \right) \setminus B'' \left(0, (1+\rho)-\frac{7}{d} \right) \\
\hline
 D'_2= d^{-1/2}(0,1) \times (0,1)
& C''_2  =  B'' \left(0, \sqrt{2}(1+\rho)-\frac{6}{d} \right) \setminus B'' \left(0, \sqrt{2}(1+\rho)-\frac{7}{d} \right)  \\
 \hline
\end{array}
$$
%\begin{eqnarray*}
%D'_0 & = & \big(-d^{-1/2},d^{-1/2}\big) \times \big( -d^{-1/2}, 0 \big), \\
%\forall i \in \{1, 2\} \quad D'_i & = & \big(0,d^{-1/2}\big)^2,
%\end{eqnarray*}
%and the followining sets in $\R^{d-2}$
%\begin{eqnarray*}
%C''_0 & = & \{0\}, \\
%\forall i \in \{1, 2\} \quad C''_i & = & B''(0, \sqrt{i}(1+\rho)-6d^{-1}) \setminus B''(0, \sqrt{i}(1+\rho)-7d^{-1}).
%\end{eqnarray*}
Finally, we set
$$
\forall i \in\{0,1,2\} \quad C_i   =   D'_i\times C''_i.
$$
 Note that for $d$ large enough, $C''_1 \cap C''_2=\varnothing$ and thus $C_1 \cap C_2=\varnothing$. The next straightforward lemma controls the asymptotics in the dimension $d$ of the volume of these sets. The proof is left to the reader.

\begin{lemma} \label{l:c}
For $i \in \{1,2\}$ :
$$
\lim_{d \to +\infty} \frac1d \ln \frac{|C''_i|}{v_{d-2}}=\lim_{d \to +\infty} \frac1d \ln \frac{|C_i|}{v_{d}}=\ln (\sqrt{i}(1+\rho)).
$$ 
\end{lemma}

%\proof \verifier{est-ce vraiment necessaire ?} Note that
%\begin{align*}
%% |C''_1| & = v_{d-2}\left( \left( (1+\rho)-\frac{6}{d} \right)^{d-2} - \left( (1+\rho)-\frac{7}{d} \right)^{d-2} \right), \\
%\frac{|C''_1|}{v_{d-2}} & = (1+\rho)^{d-2} \left( \left( 1-\frac{6}{d(1+\rho)} \right)^{d-2} - \left( 1-\frac{7}{d(1+\rho)} \right)^{d-2} \right), \\
%\frac{|C_1|}{v_{d}} & = d \frac{v_{d-2}}{v_d}\frac{|C''_1|}{v_{d-2}}=\frac{d^2}{2\pi}\frac{|C''_1|}{v_{d-2}}, 
%\end{align*}
%from which the limits for $i=1$ can be deduced. The case $i=2$ is similar.
%\finpreuve

%\medskip
We will seek the couple $(x_1,x_2)$ involved in the event ${\mathcal G}^+(x_0)$ in $C_1 \times C_2$.
But we also have to ensure that $B(x_0, \rho)\cap B(x_1,1) \ne \varnothing \text{ and } B(x_1,1)\cap B(x_{2}, \rho) \ne \varnothing$.  We set,
for $y \in C_{1}$, 
$$
D''_2(y'') =  \left\{z'' \in C''_2 : \; \left<z'',y''\right>\ge \|y''\|.\|z''\|\frac{\sqrt 2}{2}\right\} \subset \R^{d-2} \quad \text{ and } \quad 
D_2(y) =  D'_2 \times D''_2(y'')\subset C_2.
$$
The set $D''_2(y'')$ is the intersection of the annulus $C''_2$ and of a cone with axis $y''$.

\begin{lemma} \label{l:G}
1. If the dimension $d$ is large enough, 
\begin{align}
\forall y \in C_0 \quad & C_1 \subset B(y,1+\rho),  \label{e:inclus1}\\
\forall y \in C_1 \quad & D_2(y) \subset B(y,1+\rho) \cap C_2.\label{e:inclus2}
\end{align}
2. Let $x_0 \in C_0$, and take $d$ large enough to have \eqref{e:inclus1} and \eqref{e:inclus2}. If there exist $X_1 \in \chi_1 \cap C_1$ and $X_{2} \in \chi_{\rho} \cap D_2(X_1)$, then the event ${\mathcal G}_+(x_0)$ occurs.
\end{lemma}

\proof 
1. Let  $y \in C_0$ and $z \in C_1$. For $d$ large enough : 
$$
\|z-y\|^2 =  \|z'-y'\|^2+ \|z''-y''\|^2  \le \frac8d+((1+\rho)-6d^{-1})^2 < (1+\rho)^2.
$$
Let  now $y \in C_{1}$ and $z \in D_{2}(y)$. Then, as soon as $d$ is large enough,
\begin{eqnarray*}
 \|z-y\|^2 
 & = & \|z'-y'\|^2+ \|z''-y''\|^2 \\
 & \le & \frac2d+\|y''\|^2+\|z''\|^2-2<y'',z''> \\
 & \le & \frac2d+ \left(1+\rho-\frac{6}{d}\right)^2 + \left(\sqrt 2(1+\rho)-\frac{6}{d}\right)^2 -2\left(1+\rho-\frac{7}{d}\right)\left(\sqrt 2(1+\rho)-\frac{7}{d}\right) \frac{\sqrt 2}{2} \\
 & \le & (1+\rho)^2 + \frac{1}{d}\left(2 -(12-7 \sqrt2)(1+\rho)(1+\sqrt2)\right)  +O(d^{-2})  < (1+\rho)^2. 
\end{eqnarray*}

2. Let $x_0 \in D_0$. Assume there exist $X_1 \in \chi_1 \cap C_1$ and $X_{2} \in \chi_{\rho} \cap D_2(X_1)$.
Then,  
\begin{itemize}
\item[$\bullet$] \eqref{e:inclus1} ensures that $\|X_1-x_0\|<1+\rho$, and thus $B(x_0,\rho)\cap B(X_1,1)\ne \varnothing$;
\item[$\bullet$] \eqref{e:inclus2} ensures that $\|X_2-X_1\|<1+\rho$, and  thus $B(X_1,1)\cap B(X_2,\rho)\ne \varnothing$.
\end{itemize} 
Thus ${\mathcal G}_+(x_0)$ occurs. 
\finpreuve

\medskip 
The volume $|D''_2(y'')|$ does not depend on $y \in C_{1}$, and is denoted by $|D''_2|$. We now give asymptotic estimates for $|D''_2|$ :

\begin{lemma} \label{l:d}
$$
\lim_{d \to +\infty} \frac1d \ln \frac{|D''_{2}|}{v_{d-2}} = \lim_{d \to +\infty} \frac1d \ln \frac{|D_{2}|}{v_d} =\ln (1+\rho).
$$
\end{lemma}

\proof
We have, by homogeneity and isotropy:
\begin{equation} \label{e:louis}
|D''_2| = \left((\sqrt 2(1+\rho)-6d^{-1})^{d-2}-(\sqrt 2(1+\rho)-7d^{-1})^{d-2}\right) |S|
\end{equation}
where
$
S = \{x=(x_1,\dots,x_{d-2}) \in B''(0,1) : x_1 \ge \|x\| \frac{\sqrt 2}{2}\}.
$\\
But $S$ is included in the cylinder 
$$
\{(x_i)_{1 \le i \le d-2} \in \R^{d-2} : x_1 \in [0,1], \|(x_2,\dots,x_{d-2})\| \le \frac{\sqrt 2}{2}\}
$$
and $S$ contains the cone 
$$
\{ (x_i)_{1 \le i \le d-2}\in \R^{d-2} : x_1 \in [0,\frac{\sqrt 2}{2}], \|(x_2,\dots,x_{d-2})\| \le x_1\}.
$$
Therefore :
\begin{equation} \label{e:louis2}
v_{d-3} \left(\frac{\sqrt 2}{2} \right)^{d-2}(d-2)^{-1} \le |S| \le v_{d-3}\left(\frac{\sqrt 2}{2} \right)^{d-3}.
\end{equation}
>From \eqref{e:louis} and  \eqref{e:louis2}, we get 
$$ \lim_{d \to +\infty} \frac1d \ln\left(\frac{|D_2''|}{v_{d-2}}\right) = \ln(1+\rho).
$$
The lemma follows. 
Note that a direct calculus with spherical coordinates can also give the announced estimates.\finpreuve

%Estimons maintenant le volume de $D_{i}(y)$.
%Remarquons que pour tout $i \in \{1, \dots, k+1\}$, 
%$|D'_i|=\frac{v_2}{9(k+1)d}$.   
%Ainsi, 
%$\displaystyle \lim_{d \to +\infty} \frac1d \ln \left(|D'_i|\frac{v_{d-2}}{v_d}\right)=0$. 
%
%En int\'egrant en coordonn\'ees sph\'eriques, il vient
%$$|D''_i(y'')|=v_{d-2} ((d_i+d^{-1})^{d-2}-d_i^{d-2})\frac{\int_0^{\theta_i} \sin^{d-4}(\varphi)d\varphi}{\int_0^\pi \sin^{d-4}(\varphi) d\varphi}.$$
%Ainsi, (remarquons que c'est ind\'ependant de $y$)
%\begin{eqnarray*}
%& & \frac1d \ln \frac{|D''_i(y'')|}{v_{d-2}} \\
%& = & \frac{d-2}d \ln d_i + \frac1d \ln \left( (1+(dd_i)^{-1})^{d-2}-1 \right) 
%+ \frac1d \ln\int_0^{\theta_i} \sin^{d-2}(\varphi)d\varphi 
%- \frac1d \ln \int_0^\pi \sin^{d-2}(\varphi) d\varphi
%\end{eqnarray*}
%Le dernier terme (int\'egrale de Wallis, $\int_0^\pi \sin^{n}(\varphi) d\varphi \simeq \sqrt{2\pi/n}$) tend vers $0$; d'autre part,
%$$\frac1d \ln \int_0^{\theta_i} \sin^{d-2}(\varphi)d\varphi \to \ln \left(\frac{r_i\sqrt{1-a_i^2}}{d_i}\right).$$
%$\;$ \hfill $\blacksquare$
%\end{proof}

\medskip
\proofof{Proposition \ref{p:onestep}}
Choose $p<1$ and  $x_0 \in W$ such that $x''_0=0$.

$\bullet$ 
We start with a single individual, encoded by its position $\zeta_0=\{x_0\} \subset C_0$, and we set
$$
\zeta_1 = \chi_1 \cap   C_1 \quad \text{ and }\quad 
\zeta_{2} = \chi_{\rho} \cap \bigcup_{y \in \zeta_1} D_{2}(y) \subset C_{2}. 
$$ 
By Lemma \ref{l:G}, for $d$ large enough, if $\zeta_{2}\neq \varnothing$ then the event ${\mathcal G}^+(x_0)$ occurs. To bound from below the probability that $\zeta_{2} \neq \varnothing$, we build a simpler random set $\xi$, stochastically dominated by~$\zeta_2$.

$\bullet$ We set $\alpha_1=\lambda_1 |C_1|$ and $\alpha_{2}=\lambda_{\rho} |D_{2}|$ : thus, $\alpha_i$ is the mean number of children of a point in $\zeta_{i-1}$. 

Consider a random vector $X=(X_1,X_{2})$ of points in $\R^d$ defined as follows : 
$X_1$ is taken uniformly in $C_1$, then $X_2$ is taken uniformly in $D_2(X_1)$.
We think of $X$ as a potential single branch of progeny of $x_0$.
Let then $(X^j)_{j \ge 1}$ be independent copies of $X$.
Let now $N$ be an independent Poisson random variable with parameter $\alpha_1$ : this random variable gives the number $|\zeta_1|$ of children of $x_0$.
We will use the $N$ first $X^j$, one for each child of $x_0$.

We now take into account the fact that some individuals may have no children. We shall deal with geometric dependencies later. 
%Note that in our new process each individual of any generation $i, \,i \ge 1$, has at most one child. We made that choice in order to handle more easily geometric dependencies.
Let $Y=(Y^j)_{ j \ge 1}$ be an independent family of independent random variables, such that $Y^j$ follows the Bernoulli law with parameter $1-\exp(-\alpha_2)$, which is the probability that a Poisson random variable with parameter $\alpha_2$ is different from $0$. 
We set $J_1=\{1,\dots,N\}$ and
$$
J_2 = \{ 1 \le j \le N : Y^j = 1\}.
$$
Thus the random set $J_2$ gives the superscripts of the individuals, among the $N$ individuals of the first generation, that have at least one child in a process with no dependencies due to geometry.

To take into account the geometrical constraints between individuals, we set, for every  $j \ge 1$, 
\begin{align*}
& Z^j=1 \text{ if }X^j_2 \not\in \bigcup_{j' \in J_{1} \setminus \{j\}} D_2(X_{1}^{j'}) \quad \quad \text{ and } Z^j=0 \text{ otherwise}, \\
&\xi=\{X_2^j :  j \in J_2 \hbox{ and } Z^j=1\}.
\end{align*}
We thus reject an individual $X_2^j$ as soon as $Z^j=0$.
Recall that, when building generation $2$ from generation $1$, we explore the Poisson point processes in the area $\bigcup_{j \in J_{1}} D_2(X_{1}^{j})\subset C_2$. 
Remember that by construction, $C_1$ and $C_2$ are disjoint. 
Therefore, one can check that the set
$\xi$
is stochastically dominated  by $\zeta_2$\footnote{ Note that the random set $\{X_1^1,\dots,X_1^N\}$ has the same distribution as $\zeta_1$. In order to build a random set with the same distribution as $\zeta_2$, we could proceed as follows.  Let $(N^j)_{j \ge 1}$ be independent random variable distributed according to the Poisson distribution with mean $\alpha_2$. Throw $N^1$ random points uniformly in $D_2(X_1^1)$. Then throw $N^2$ random points uniformly in $D_2(X_1^2)$ and remove the points that fell in $D_2(X_1^1)$. Then throw $N^2$ random points uniformly in $D_2(X_1^3)$ and remove the points that fell in $D_2(X_1^1)$ or in $D_2(X_1^2)$. And so on. The random set of all the points thrown and not removed has the same distribution as $\zeta_2$. 

In the proof of Proposition 2.3, we reject more points than in this classical construction, thus only obtaining a stochastic domination: 
\begin{itemize}
\item First, we replace $N^j$ by $\min(1,N^j)$ to keep at most one point $X_2^j$ for each $j$ (this is the role of $Y^j$).
\item Secondly, we reject this point $X_2^j$ as soon as it falls into any of the $D_2(X_1^{j'})$ for $j'\ne j$ instead of only forbiding the $D_2(X_1^{j'})$ for $j'<j$ (this is the role of $Z^j$).
\end{itemize}
}. Thus to prove Proposition \ref{p:onestep}, we now need to bound from below the probability that $\xi$ is not empty.

$\bullet$ Let $T$ be the smallest integer $j$ such that $ Y^j=1$ : in other words, $T$ is the smallest superscript  of a branch that lives till generation $2$. To ensure that $\xi\neq \varnothing$, it is sufficient that $T\le N$ and that $Z^T =1$. So :
\begin{eqnarray*}
1-P({\mathcal G}^+(x_0)) 
 & \le & P(\xi=\emptyset) 
  \le  P(\#J_{2} = 0) + P \left(\{T \le N\} \cap   \{Z^T = 0\}\right) .
\end{eqnarray*}
By construction:
\begin{eqnarray*}
P(T \le N \hbox{ and } Z^T=0) 
 & = & P\left(T \le N, \; \exists j \in J_{1} \setminus \{T\}  \hbox{ such that } X_2^T \in D_2(X_{1}^j) \right) \\
 & \le & \sum_{j \ge 1} P\left(T \le N \hbox{ and } j \in J_{1} \setminus \{T\}  \hbox{ and } X_2^T \in D_2(X_{1}^j) \right) \\
 & = & \sum_{j \ge 1} E\left( 1_{T \le N} 1_{j \in J_{1} \setminus \{T\}} P\left( X_2^T \in D_2(X_{1}^j) \left.\right| Y,N \right)\right) \\
 & = & \sum_{j \ge 1} E\left( 1_{T \le N} 1_{j \in J_{1} \setminus \{T\}} \right)P\left( X_2^1 \in D_2(X_{1}^2) \right) \\
 & \le & E(\# J_{1})P\left( X_2^1 \in D_2(X_{1}^2) \right)=E(N)P\left( X_2^1 \in D_2(X_{1}^2) \right).
\end{eqnarray*}
Besides, as  $(X_2^1)''$ is uniformly distributed on  $C''_2$  and is independent of $(X_{1}^2)''$, 
\begin{equation}
P\left(X_2^1 \in D_2(X_{1}^2)\right) 
  =   P\left((X_2^1)'' \in D''_2((X_{1}^2)'')\right) 
  =  \frac{|D''_2|}{|C''_2|} \label{e:probainterference} . \nonumber
\end{equation}
This leads to 
\begin{equation}
1-P({\mathcal G}^+(x_0)) \le P(\#J_{2} = 0) +  E(N)\frac{|D''_2|}{|C''_2|}. \label{e:majJ}
\end{equation}

$\bullet$ $N$ follows a Poisson law with parameter $\alpha_1=\lambda_1 |C_1|$ with $\lambda_1=\frac{\kappa^d}{v_d 2^d}$. Thus
$$E(N)=\frac{\kappa^d}{2^d} \frac{|C_1|}{v_d}.$$
Lemmas \ref{l:c} and \ref{l:d} ensure that :
$$
\lim_{d\to +\infty} \frac1d  \ln \left(\frac{|C_1|}{v_d}\right) = \ln(1+\rho) \quad \text{ and } \quad \lim_{d\to +\infty} \frac1d  \ln \left(\frac{|D''_2|}{|C''_2|}\right) = \ln\left(\frac{1}{\sqrt 2}\right).
$$
Thus,  we have :
\begin{eqnarray}
&&\lim_{d\to +\infty} \frac1d \ln \left( E(N)\frac{|D''_2|}{|C''_2|}\right) 
\le  \ln \left( \frac {(1+\rho) \kappa}{2\sqrt 2}\right)<0 \quad \text{since }\kappa<\frac{2 \sqrt 2}{1+\rho}, \nonumber \\
\text{therefore } &&\lim_{d\to +\infty}  E(N)\frac{|D''_2|}{|C''_2|}=0.
\label{e:deuxiemeterme}
\end{eqnarray}
The cardinality of $J_2$ follows a Poisson law with parameter
$$
\eta = \alpha_1  (1-\exp(-\alpha_{2})).
$$
Remember that  $\alpha_1=\lambda_1 |C_1|$, $\alpha_{2}=\lambda_{\rho} |D_{2}|$, $\lambda_1=\frac{\kappa^d}{v_d 2^d}$ and $\lambda_{\rho}=\frac{\kappa^d}{v_d 2^d \rho^d}$. By Lemma \ref{l:d}, we have the following limits:
\begin{eqnarray*}
\lim_{d \to +\infty} \frac1d \ln \alpha_1 & = & \ln \frac{\kappa(1+\rho)}{2} >  0, \\
\lim_{d \to +\infty} \frac1d \ln  \alpha_{2} & = & \ln(\kappa\frac{1+\rho}{2\rho}) < 0.
\end{eqnarray*}
The first inequality is a consequence of $\kappa > \kappa_\rho^c$.
The second inequality is a consequence of $\kappa < 1$.
Consequently, we first see that 
\begin{eqnarray}
\lim_{d\to +\infty}  \frac1d \ln(\eta) 
 & = & \lim_{d\to +\infty}  \frac1d \ln(\alpha_1\alpha_{2})  
  =  \ln \left(\kappa^{2}\frac{(1+\rho)^2}{4\rho} \right)> 0  \nonumber \\
 \text{therefore,} &&\lim_{d\to +\infty} P(\#J_{2}=0) = 0. \label{e:premierterme}
\end{eqnarray}
The inequality is a consequence of $\kappa > \kappa_\rho^c$.

To end the proof, we put estimates \eqref{e:premierterme} and  \eqref{e:deuxiemeterme} in \eqref{e:majJ}.
\finpreuve

\subsubsection{Several steps in the $2$-dimensional oriented percolation model}
\label{s:severalsteps}

We prove here Proposition \ref{p:supercritical} by building the supercritical $2$-dimensional oriented percolation process embedded in the two-type Boolean Model.

\proofof{Proposition \ref{p:supercritical}}
We first define an oriented graph in the following manner: the set of sites is 
$$S=\{(a,n)\in \Z \times \N: \; |a| \le n, \; a+n \hbox{ is even }\};$$
from any point $(a,n) \in S$, we put an oriented edge to $(a+1,n+1)$, and an oriented edge to $(a-1, n+1)$. We denote by 
$\vec{p}_c \in (0,1)$ the critical parameter for Bernoulli percolation on this oriented graph -- see Durrett \cite{Durrett-oriented} 
 for results on oriented percolation in dimension~2.  

For any $(a,n) \in S$, we define the following subsets of $\R^d$ %\verifier{justifier le $d^{-1/2}$} :
\begin{eqnarray*}
W_{a,n} & = &  d^{-1/2} \left((a-1,a+1) \times (n-1,n) \times \R^{d-2}\right), \\
W^-_{a,n} & = &  d^{-1/2} \left((a-1,a) \times (n,n+1) \times \R^{d-2}\right), \\
W^+_{a,n} & = &  d^{-1/2} \left((a,a+1) \times (n,n+1) \times \R^{d-2}\right).
\end{eqnarray*}
Note that the $(W_{a,n})_{(a,n) \in S}$ are disjoint and that $W^+_{a,n}\cup W^-_{a+2,n} \subset W_{a+1,n+1}$.

We now fix  $\kappa \in (\kappa^c_{\rho},1)$, and 
for $x_0 \in W_{a,n}$, we introduce the events :
\begin{eqnarray*}
{\mathcal G}^+_{a,n}(x_0) & = & \left\{
\begin{array}{c}
\hbox{There exist } x_1 \in \chi_1 \cap  W^+_{a,n} \hbox{ and } x_{2} \in \chi_{\rho} \cap W^+_{a,n} \\
\hbox{such that } B(x_0, \rho)\cap B(x_1,1) \ne \varnothing \text{ and } B(x_1,1)\cap B(x_{2}, \rho) \ne \varnothing
\end{array}
\right\}, \\
{\mathcal G}^-_{a,n}(x_0) & = & \left\{
\begin{array}{c}
\hbox{There exist } x_1 \in \chi_1 \cap  W^-_{a,n} \hbox{ and } x_{2} \in \chi_{\rho} \cap W^-_{a,n} \\
\hbox{such that } B(x_0, \rho)\cap B(x_1,1) \ne \varnothing \text{ and } B(x_1,1)\cap B(x_{2}, \rho) \ne \varnothing\end{array}
\right\}.
\end{eqnarray*}
Note that ${\mathcal G}^+_{0,0}(x_0)$ is exactly the event ${\mathcal G}^+(x_0)$ introduced in \eqref{e:bonevnt}, and that the other events are obtained from this one by symmetry and/or translation.

Next we choose $p \in (\vec{p}_c, 1)$. With Proposition \ref{p:onestep}, and by translation and symmetry invariance, we know that for every large enough dimension $d$, for every $(a,n) \in S$, for every $x \in W_{a,n}$:
\begin{equation}\label{e:gp}
P({\mathcal G}^{\pm}_{a,n}(x))\ge p.
\end{equation}
We fix then a dimension $d$ large enough to satisfy \eqref{e:gp}.
We can now construct the random states, open or closed, of the edges of our oriented graph. The aim is to build inductively some appropriate paths of balls from a ball centered at a point $x(0,0) \in W_{0,0}$ to balls centered at points $x(a,n) \in W_{a,n}$. In case of failure for a given $(a,n)$, we find it convenient to set $x(a,n)=\infty$, where $\infty$ denotes a virtual site. In the end, usefull paths will only use finite $x(a,n)$.
%\verifier{We denote by . This virtual site is only introduced to build an explicite coupling between our dependent oriented percolation and independent oriented percolation. We can check at the end of the construction that an open path starting from $(0,0)$ never visits this virtual site.}

\prg{Definition of the site on level $0$}
Almost surely, $\chi_\rho \cap W_{0,0} \neq \varnothing$. 
We take then some $x(0,0) \in \chi_\rho \cap W_{0,0}$. 

\prg{Definition of the edges between levels $n$ and $n+1$} Fix $n \ge 0$ and assume we have built a site $x(a,n) \in W_{a,n}\cup \{\infty\}$  for every $a$ such that $(a,n) \in S$.  Consider $(a,n) \in S$ :
\begin{itemize}
\item If $x(a,n)=\infty$ :  
we decide that each of the two edges starting from $(a,n)$ is open with probability $p$ and closed with probability $1-p$, independently of everything else; we set $z^-(a,n)=z^+(a,n)=\infty$.
\item Otherwise, $x(a,n) \in W_{a,n}$ and : 
\begin{itemize}
\item Edge to the left-hand side :
\begin{itemize}
\item if the event ${\mathcal G}^-_{a,n}(x(a,n))$ occurs : 
we take for $z^-(a,n)$ some point $x_2\in {W}^-_{a,n} \subset {W}_{a-1,n+1}$ given by the occurrence of the event, and we open the edge from $(a,n)$ to $(a-1,n+1)$  ; 
\item otherwise : we set $z^-(a,n)=\infty$ and we close the edge from $(a,n)$ to $(a-1,n+1)$. 
\end{itemize}
\item Edge to the right-hand side : 
\begin{itemize}
\item if the event ${\mathcal G}^+_{a,n}(x(a,n))$ occurs : 
we take for $z^+(a,n)$ some point $x_2\in {W}^+_{a,n} \subset {W}_{a+1,n+1}$ given by the occurrence of the event, and we open the edge from $(a,n)$ to $(a+1,n+1)$  ; 
\item otherwise : we set $z^+(a,n)=\infty$ and we close the edge from $(a,n)$ to $(a+1,n+1)$. 
\end{itemize}
\end{itemize}
\end{itemize}
For $(a,n)$ outside $S$, we set $z^{\pm}(a,n)=\infty$.

\prg{Definition of the sites at level $n+1$} Fix $n \ge 0$ and assume we determined the state of every edge between levels $n$ and $n+1$. Consider $(a,n+1) \in S$ :
\begin{itemize}
\item If $z^+(a-1,n) \neq \infty$ : set $x(a,n+1)=z^+(a-1,n)\in {W}_{a,n+1}$.
\item Otherwise :
\begin{itemize}
\item if $z^-(a+1,n) \neq \infty$ : set $x(a,n+1)=z^-(a+1,n) \in {W}_{a,n+1}$,
\item otherwise : set $x(a,n+1)=\infty$.
\end{itemize}
\end{itemize}

\bigskip
Assume that there exists an open path of length $n$ starting from the origin in this oriented percolation : we can check that the leftmost open path of length $n$ starting from the origin gives a path in the two-type Boolean model along which   balls with radius $1$  alternate with  balls with radius $\rho$. Thus, percolation in this oriented percolation model implies alternating percolation in the two-type Boolean model. Let us check that percolation occurs indeed with positive probability.

For every $n$, denote by ${\mathcal F}_n$ the $\sigma$-field generated by the restrictions of the Poisson point processes $\chi_1$ and $\chi_{\rho}$ to the set
$$
d^{-1/2} \left( \R \times (-\infty,n) \times \R^{d-2} \right).
$$
By definition of the events $\mathcal G$ -- remember that the  $(W_{a,n})_{(a,n) \in S}$ are disjoint -- and by \eqref{e:gp}, the states of the different edges between levels $n$ and $n+1$ are independent conditionally to ${\mathcal F}_n$. Moreover, conditionally to ${\mathcal F}_n$, each edge between levels $n$ and $n+1$ has a probability at least $p$ to be open. 
Therefore, the oriented percolation model we built stochastically dominates Bernoulli oriented percolation with parameter $p$. As $p>\vec{p}_c$, with positive probability, there exists an infinite open path in the oriented percolation model we built; this ends the proof of Proposition \ref{p:supercritical}.
\finpreuve

\subsection{Proof of Theorem \ref{t:2}}
If $\kappa<\kappa^c_\rho$ then, by Proposition \ref{p:subcritical}, there is no percolation for $d$ large enough.
Therefore, for any such $\kappa$ and for any large enough $d$ we have:
$$
\lambda^c_d(\mu_d) \ge \frac{\kappa^d}{v_d 2^d} \quad \text{and then} \quad \widetilde{\lambda}^c_d(\mu_d) = \lambda_d^c(\mu_d) v_d 2^d\int r^d \mu_d(dr)\ge 2\kappa^d .
$$
Letting $d$ go to $+\infty$ and then $\kappa$ go to $\kappa^c_\rho$, we then obtain
\begin{equation}
\label{e:liminf1}
\liminf_{d \to +\infty}\frac1d \ln \left(\lambda^c_d(\mu_d)\right)  \ge \ln \left( \kappa^c_\rho\right) .
\end{equation}

Choose now $\kappa$ such that $\kappa^c_\rho<\kappa$. 
By Proposition \ref{p:supercritical}, there is percolation for $d$ large enough in $\Sigma$.
Therefore, for any $\kappa>\kappa^c_\rho$ and for any large enough $d$ we have, as before:
$$
\lambda^c_d(\mu_d) \le \frac{\kappa^d}{v_d 2^d} \quad \text{and then} \quad\widetilde{\lambda}^c_d(\mu_d) \le 2\kappa^d.
$$
Letting $d$ go to $+\infty$ and then $\kappa$ go to $\kappa^c_\rho$, we then obtain
\begin{equation}
\label{e:limsup1}
\limsup_{d \to +\infty}\frac1d \ln \left(\lambda^c_d(\mu_d)\right)  \le \ln \left( \kappa^c_\rho\right) .
\end{equation}
Bringing \eqref{e:liminf1} and \eqref{e:limsup1} together, we end the proof of Theorem \ref{t:2}. \finpreuve

\def\cprime{$'$} \def\cprime{$'$}

\end{document}